# Probabilistic Events and Physical Reality:
# a Complete Algebra of Probability

Paolo Rocchi and Leonida Gianfagna
*IBM, via Shangai 53, 00144 Roma, Italy*

Resumé: Ce travail dérive d'un etude plus large concernant les fondements de la probabilité qui a la caractéristique suivant. Nous n'abordons pas tout de suite la probabilité, comme il est d'usage parmi les auteurs, mais nous commençons par l'analyse de l'argument de la probabilité qui est l'événement probabiliste.

Nous débutons par une discussion qui regarde deux modéles formales: la phrase linguistique et l'ensemble qui sont le représentations les plus importantes de les événements physique que l'on observe dans le monde. Tous les deux ont de contradictions et de cela nous déduons que la rigoreuse définition de l'événement c'est la question basilaire et inévitable des fondements.

Nous proposons la structure des niveaux pour modeler l'événement indéterministe. Cette structure algébrique suggère une définition de la probabilité qui a une explicite signifié physique. De plus elle harmonise les interpretations subjectives et objectives de la probabilité qui l'une avec l'autre sont en opposition jusqu'ici.

Abstract: This contribution derives from a rather extensive study on the foundations of probability. We start by discussing critically the two main models of the random event in Probability Theory and cast light over a number of incongruities. We conclude that the argument of probability is the critical knot of the probability foundations and put forward the structure of levels for the partially determinate event. The structural model enables us to define the probability and to attune its subjective and objective interpretations.

Keywords: Probability foundations, Random events.

## I. FOREWORD

The Probability Calculus dates back to Pascal who in about 1654 began to investigate the chances of getting different values for rolls of the dice. The talented Frenchman did not close his mind to particular problems but addressed himself to a much wider and complete scenario. Pascal intended to place all probabilistic themes into the ad hoc sector that he outlined as follows: "By joining the rigor of the scientific demonstrations with the uncertainties of chance and by reconciling these two apparently contradictory items, we can assign the following amazing title: *The Geometry of Chance* to the discipline embracing both of them."

Since then scientists conceived random events as "points" in this fascinating and non-obvious geometry. Several theories were proposed for probability and each of them confirms the event as the argument of probability but few scholars (to the best of our knowledge) are explicitly interested in discussing it[1,2]. Authors debate fiercely the meaning of probability but its argument is taken as a secondary matter. Somebody even explicitly declares the event is irrelevant from the mathematical point of view [3].

A logical order and an obvious rule require that at first we must discuss the argument and then its measure. In particular the open debate on the probability interpretation would require this rigid method. Instead few specialists check out Pascal's conjecture, they rarely justify the theoretical models of the random event in relation to the experimental data, they do not discuss critically the formal representations in this or that practical circumstance, when they conform more or less with the physical experiment.

We shall illustrate our remarks on the theoretical models of the event, which are substantially two. Kolmogorov and many authors conceive the event as a subset. The logical and the subjective schools, more sensitive to the experienced disparities, adopt the linguistic model for the events. We shall highlight some incongruities which solicit us to go deep into the argument of probability. We shall put forward a solution to the problem.

## II. LINGUISTIC MODEL

Let the propositions *p, q, ...* describe the same random event, thus they are equivalent

$$p \Leftrightarrow q \qquad (1)$$

and form the equivalence class *X*

$$X = \{p, q, ...\} \qquad (2)$$

that constitutes the argument of the probability

$$P = P(X) \qquad (3)$$

Various theories assume special formalization reflecting the algebraic, the logic etc. formulations [4,5]. They do not interest us, instead we focus on the linguistic model (2).

**Remark 2.1:** Some subjectivist theorists acutely hold that any event is intricate and the linguistic model is consistent with this well experienced complication. Disciplines treating complex phenomena such as psychology and sociology, business management and medicine adopt the linguistic representation and consider other schemes to be reductive. Therefore, the proposition seems an adequate model except for the following perplexity: each primitive is a simple idea and can be left to intuition only for this reason. For example a *number,* a *point,* an *entity* are elementary notions. Can we declare that the random event is complicated and simultaneously assume it as a primitive ?

Complexity opposes the assumption of a primitive. This contrast would at least require an in depth justification that instead is lacking, as far as we know.

**Remark 2.2:** For decades linguistic specialists claim the sentences *p, q, ...* written in natural language may be equivocal. Some words are ambiguous to the extent that they must be interpreted. Bayesians exclude this risk and assume that a proposition either represents one event or does not represent it. We prove this assumption is unwarranted by means of Bertrand's paradox. The following proposition

*p* = "A chord is longer than the side of the equilateral triangle inscribed within the circle"  (4)

portrays a physical event and according to the Bayesan method we write

$$p = \text{true} \qquad (5)$$

The relation (5) holds that *p represents one physical event* but the description *p* conceals *six* different dynamics and we shall describe two of the best known.

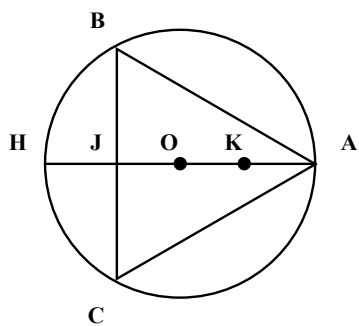

Fig. 1

**2.2.1:** Let the equidistant points H, J, O, K and A in the diameter and the chord fall parallel to BC. Under hypothesis (4), the chord proves to be greater than the side of the triangle when it falls between J and K, while the chord is less than the side when it falls in HJ and in KA. Since any segment is equally probable, the probability that the cord is longer than the triangle side is 1/2

**2.2.2:** Now we fix an extreme of the chord in A, while the other extreme H moves on the circumference. The triangle ABC is equilateral, the arcs AB, BC and CA are equal and the chord has the same probability of having its extreme H in one of the three arcs. Only if H falls on the arc BC then the chord is greater than the side of the equilateral triangle, therefore the probability is 1/3

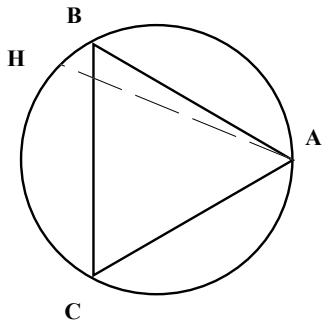

Fig. 2

The sentence (4) does not picture a sole physical event and yields two values of probability. We conclude that in general the linguistic expression is equivocal and constitutes a problematic model for the probability calculation.

## III. SET MODEL

Kolmogorov assumes that the sample space $\Omega$ includes all the possible elementary events [6]. He defines the random event $X$ as a set of specific elementary events

$$X = \{\xi_1, \xi_2, \xi_3, ...\} \qquad (6)$$

where $X$ is a subset of $\Omega$

$$X \subset \Omega \qquad (7)$$

and the probability measures $X$
.
$$P = P(X) \qquad (8)$$

**Remark 3.1:** Kolmogorov immediately clarifies the application of his theory. In Section 2 of Chapter 1 he interprets $X$ as "the set of all possible variants $\xi_1, \xi_2, \xi_3,...$ of the outcome of the given event". Countless scholars share the same conception, which causes some perplexities in the light of the modern systemic knowledge[7]. Any physical event is a dynamic evolution beginning with its antecedent and ending with its consequent. E.g. Drawing a card constitutes an event whose antecedent is the whole pack of cards, and the consequent is the selected card. The result, namely the consequent, is a part and the event is the whole

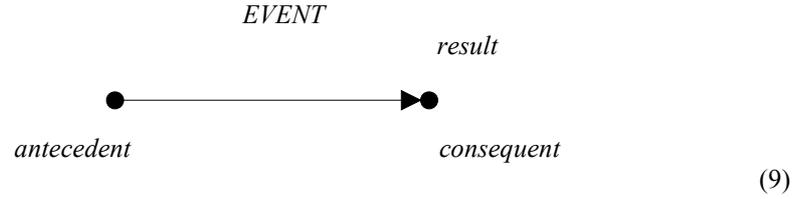

$$\text{EVENT}$$
antecedent — result → consequent (9)

The properties of the event quite differ from the properties of the outcome and we find an open inconsistency. We cannot call $X$ "set of events " and also conceive it as a "set of results". We cannot merge the part and the whole without a logical justification but this is lacking so far.

**Remark 3.2:** Definition (8) entails that a result emerges in any application. Aside from gambling and some other special events, this statement encounters difficulties because outcomes may be lacking. E.g. Mr. Tom gets information about Mr. John. He takes a decision, which influences his inner psychological state, but does not tell his resolution. Nothing emerges outside his mind and the event does not output anything. In this situation $X$ cannot be defined and the Kolmogorov theory cannot be applied.

**Remark 3.3:** When the result appears, sometimes it cannot be modeled as a set. Quantum Physics brings up an exception of significant importance because the outcomes in the "two slits experiments" cannot be calculated as sets [8].

## IV. DENOTED EVENTS

Some may object that the above-mentioned difficulties do not exist in the Probability Calculus or at least not to the extent we have just highlighted. They argue that the random event $X$ is the argument of $P$, thus it is conceived as a "unicum" and $X$ just denotes the practical event; $X$ codes the facts simply. This assumption justifies the generic proposition and the set referring to the result. The model $X$ serves only to identify the argument of the probability and our remarks would be out of place. Objectivists and subjectivists, bayesians, frequentists, logicists and still others amply share these considerations, therefore we plunge into the topic. We wonder: Did we verify if the random event is the rigorous argument of the probability?

Definitions (3) and (8) claim that probability depends on the event $X$, however the notion of "argument" is somewhat elaborate. We shall analyze it by discussing the speed $v$ of the moving body $x$.

**Remark 4.1:** This physical parameter depends on the space $\Delta s$ covered by $x$ during the interval $\Delta t$. The rigorous definition states the speed equals the ratio of the two values

$$v = v(s,t) = \frac{\Delta s}{\Delta t} \qquad (10)$$

**Remark 4.2:** If we detect the speeds of some bodies by an instrument, we identify the speed using the argument $x$.

$$v = v_x \qquad (11)$$

This expression points out that $v$ refers to $x$ and not to $y$ or $z$. This dependence does not make explicit the meaning of the speed but simply relates the movement to the body $x$.

We conclude that the rigorous argument in (10) provides the substantial significance of $v$ whereas (16) merely denotes the running body. The specific argument of a measure explicates the cause-effect occurring in the world. It provides the correct value through calculus whereas the generic argument does not. The Probability Calculus supplies the probability through the number of possible cases $n_p$ and the favorable cases $n_f$

$$P = \frac{n_p}{n_f} \qquad (12)$$

It gives results through the probabilities of partial events and random variables. We never calculate the probability by the sentence $X$ neither the set $X$, hence (3) and (8) are generic and do not establish the precise argument of the probability. This statement brings us to severe conclusions: mathematicians have to discover the argument of the probability before the probability itself.

In the next section, we shall discuss the solution, which we suggest for the above problems.

## V. STRUCTURAL MODEL

Ludwing von Bertalanffy, father of the General Systems Theory, conceives a system, and consequently an event, as an intricate set of items, which affect one another. Interacting and connecting is the essential character and the inner nature of events, and we take this idea at the base of our theoretical proposal.

**Postulate 5.1:** The idea of relating, of connecting, of linking is primitive.
This idea suggests two elements specialized in relating and in being related

**Definition 5.2:** The relationship $R$ connects the entities and we say $R$ has the property of connecting.

**Definition 5.3:** The entity $E$ is connected by $R$ and we say $E$ has the property of being connected.

Relationships and entities are already known as operations and elements in Algebra; as arrows and objects; as edges and vertices. The main difference is that all of them are given as primitive[9] while $R$ and $E$ derive from the axiomatic concept 5.1. In other words, the properties of the algebraic elements are openly given in 5.2 and 5.3, while they are implicit in other theories.

**Definition 5.4:** The relationship $R$ links the entity $E$ and they give the algebraic structure

$$S = (E; R) \qquad (13)$$

Note this set has the associative/dissociative property, namely the structure is the unicum $S$. Then it is defined in terms of the details $E$ and $R$. If this analysis is insufficient we reveal the entities $(E_1, E_2..., E_m)$ and the relations $(R_1, R_2, ... R_p)$; these are exploded at a greater level, and so forth.

$$\begin{array}{ll} \text{level } 0 & S \\ \text{level } 1 & E; R \\ \ldots & \ldots\ldots \\ \text{level } (n-1) & E_1, E_2..., E_m; R_1, R_2, \ldots R_p \\ \text{level } n & E_{12}....., E_{m1}, E_{m2}.., E_{mk}; R_{11}, R_{12}....., R_{p1}, R_{p2}.., R_{ph} \end{array} \qquad (14)$$

**Definition 5.5:** The structure of levels consists of the finite number $n$ of levels which are indexed and ordered sets. The element $X_j$ of the level $j$ has this property

$$X_j = \{Y_{j+1}, Z_{j+1}, W_{j+}, \ldots\} \qquad j > n \qquad (15)$$

Where $Y_{j+1}, Z_{j+1}, W_{j+1}$ are generic elements of the level $j + 1$.

In the physical reality, an event is a dynamic phenomenon linking $E_{in}$ to $E_{out}$ and from (14) we can deduce this very common structure

$$S = (E; R) = (E_{in}, E_{out}; R) \qquad (16)$$

Using a graph we get

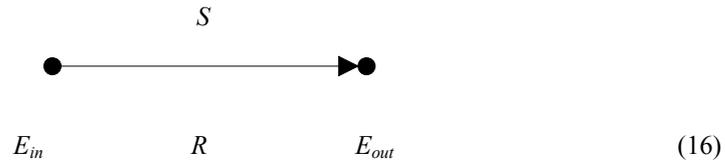

$$E_{in} \qquad R \qquad E_{out} \qquad (16)$$

The structural model identifies accurately each physical item; in particular we highlight the ensuing advantages.

**Remark 5.6:** The parts and the whole are logically separate namely the result $E_{out}$ is distinct from the event $S$ and they give a precise answer to objection 3.1. The structure can model an event, which does not produce any output in the reality and resolve question 3.2.

**Remark 5.7:** The outcome may be a rational or an irrational number, a real or an imaginary value notably $E_{out}$ can be calculated by a wave function or by a series etc. and this offers a formal solution to point 3.3. If the result $E_{out}$ is an ensemble of variables

$$E_{out} = \{\xi_1, \xi_2, \xi_3, \ldots\} \qquad (17)$$

Than the structure includes the Kolmogorov model (6).

**Remark 5.8:** A simple sentence includes nouns that are entities and a verb representing a dynamical evolution. E.g. The entities and relationship detail the following proposition

$$\text{"The coin | comes down | heads"}$$
$$E_{in} \qquad R \qquad E_{out} \qquad (18)$$

In short, the algebraic structure accomplishes the linguistic model. However a sentence can be equivocal whereas the structure $S$ is a rigorous formalism and answers remark 2.2.

**Remark 5.9:** The multiple level decomposition is amply introduced in theoretical and applied works [10]. The progressive expansion of the event is already known in the Probability Calculus where we use trees connecting the parts and the subparts of a random event [11]. Instead, the structure of levels does not require individual connections between the elements, which are inessential to the Probability Calculus.

**Remark 5.10:** Probability is the answer to questions such as: Who will win the next football game? Who will be voted in the regional elections? Shall I pass the examinations? Where is the photon now ?

These questions prove that probability is concerned with the particulars of an event that is already known in the whole. We see the overall random phenomenon but, however, we ignore the details that will produce the result. When we ask "who will win the next game?", we are familiar with the match, we already know the teams that will play, where the match will be held, etc. We comprehend the overall event, however, we do not have the details that will determine the result. Our understanding is incomplete and we translate this notion into the formalism just introduced. Let the event $S$ have the level 1, the level 2, up to the level $n$; two cases arise now.

**5.10.1:** The event running in the physical reality is entirely described by the relations and the entities in (18). The elements at level $(n+1)$ are not existent both in the model and in nature. *This structure is certain* because it wholly defines the facts.

**5.10.2:** Let the event running in the physical reality be not entirely described by the relations and the entities of the level $n$. The microelements of the level $(n+1)$ exist and influence the final results in a decisive way, however (14) cannot exhibit them. We call *uncertain* such *incomplete structure*. E.g. The structure of throwing a dice lists the dice $E_m$, the launching/falling dynamics $R_m$. The results $E_1, E_2, E_3, \ldots E_6$ and the relations, which are alternative and produce them in general, appear at the next level

$$\begin{array}{ll} \text{level } 0 & S \\ \text{level } 1 & E_m; R_m \\ \text{level } 2 & E_1, E_2, \ldots E_6\,;\, R_1, R_2, \ldots R_6, OR \\ \text{level } 3 & \ldots\ldots \end{array} \qquad (19)$$

The subrelationships of level 3 which intervene in the physical reality and cause any specific outcome are essential. They would enable the calculation of any trial and should be listed, unfortunately we are unaware of them and (19) is an uncertain structure.

**Remark 5.11:** We ignore the details of the uncertain structure and we observe that the physical event does not follow any procedure, thus we shall say it is a *random event*.

## VI. PROBABILITY AND ITS PHYSICAL EVIDENCE

The structure (19) proves how any event is equipped with precise macromechanisms and micromechanisms. The relationship $R_j$ at level $j$ runs through the subrelationships of level $(j+1)$ and any event like an industrial apparatus, a mechanical clock or an electronic device including various working parts. A certain event is entirely explained through the levels since the structure (19) clearly shows "how" the event runs through $n$ levels, which are exhaustive by definition. On the contrary, the uncertain structure is incomplete and we cannot describe "how" the event runs in the physical reality. As the impossibility of an analytical function capable of describing "how" the physical event functions since we miss the $(n+1)$ elements, we search "when" the uncertain event behaves, namely how many times it occurs in the world.

The structure of levels (16) highlights how $R$ is the pivotal element, while $E_{in}$, $E_{out}$ may be absent. The event $S$ occurs if and only if the relationship $R$ works. We measure its ability to connect that directly reflects the occurrence of $S$.

**Definition 6.1:** When $R$ always links the input to the output in the physical reality, the event $S$ is certain and the measure $P(R)$ equals one

$$P(R) = 1 \qquad (20)$$

When $R$ does not work in the physical reality, $S$ is impossible in the facts and the measure $P(R)$ is zero

$$P(R) = 0 \qquad (21)$$

If $R$ occasionally links $E_{in}$ to $E_{out}$, the connection is neither sure nor impossible and $P(R)$ assumes a decimal value

$$0 < P(R) < 1 \qquad (22)$$

We call *probability* the measure $P(R)$ of the operation $R$ which extensively indicates the occurrence of $S$. It quantifies the possibility or the impossibility of the uncertain event and meets with common sense. E.g. We throw a dice, the relationship $R_2$ sometimes connects the input to the output $E_3$. This operation occasionally works, the event is uncertain and $P(R_3)$ is decimal. E.g. The gravitation force $R_g$ attracts constantly the Earth $E_E$ to the Sun $E_S$. This relationship operates always, the event is certain and the probability $P(R_S)$ is unit.

**Remark 6.2:** We can denote the operation by means of its outcome in some special events. If we assume this univocal relation between $E_{out}$ and $R$

$$E_{out} \Rightarrow R \qquad (23)$$

Than we calculate the probability of the outcome

$$P(E_{out}) = P(R) \qquad (24)$$

E.g. The result heads $E_h$ appears whenever $R_h$ works, the output entails the operation

$$E_h \Rightarrow R_h \qquad (25)$$

And we calculate the probability using the outcome as a variable

$$P(E_h) = P(R_h) = 0.5 \qquad (26)$$

In conclusion if (17), (24) and (25) are true, Definition 6.1 consists with Kolmogorov's theory and evidences how this theory cannot calculate all the possible cases of the physical reality.

**Remark 6.3:** The accurate modeling of a physical event by means of a structure avoids ambiguities and misunderstandings in the probability calculus. For ease, Bertrand's paradox lists two different dynamics in 2.2. The chord $E_x$, antecedent to the fall, and the chord dropped in the circle $E_y$ are identical in 2.2.1 and 2.2.2. The two physical dynamics entail two relationships and then two different structures

$$\begin{aligned} S_1 &= (E_x, E_y; R_1) \\ S_2 &= (E_x, E_y; R_2) \end{aligned} \qquad (27)$$

The falls $R_1$ are parallel to the diameter, $R_2$ tie an extreme of the chord in A, thus

$$P(R_1) \neq P(R_2) \qquad (28)$$

**Remark 6.4:** Certain structures include only certain elements, impossible elements have no sense and are omitted. The unitary value of probability merely confirms what is already related in the levels. Conversely the uncertain structure lacks the lowest elements that are essential and the probability provide information about them. The decimal values of $P(R)$ clarify the intervention of the elements at level $(j + 1)$ although they are absent in the structure. For example, we ignore the parts of $R_h$ producing the result $E_h$, but (26) explains how they work. Exactly half of the $S$ occurrences is due to the subrelationships of tails $R_t$ and the other half is activated by the components of heads $R_h$. The explicatory and predictive values of probability appear absolutely relevant.

**Remark 6.5:** Definition 6.1 forces us to prove its physical meaning since it is suggested by experience. In order to simplify the discussion let the uncertain event include either the relationship $R_i$ or $NOT\ R_i$ at level $n$

$$\begin{aligned} &\text{level } 0 &&S \\ &\text{level } .. &&E; R \\ &\text{level } n &&E_i, NOT\ E_i; (R_i,\ NOT\ R_i) \\ &\text{level} &&........ \end{aligned} \qquad (29)$$

The probability $P(R_i)$ expresses the runs of $R_i$ by definition, thus the *occurrences* $gs(R_i)$ in the sample $s$ is the experimental measurement corresponding to $P(R_i)$. As much as $R_i$ connects, so much is $gs(R_i)$. Vice versa as little $R_i$ runs, so small is $gs(R_i)$. As the *absolute frequency* $gs(R_i)$ exceeds the range [0,1], we select the *relative frequency* $Fs(R_i)$ which verifies

$$0 \leq Fs(R_i) \leq 1 \qquad (30)$$

The relative frequency must coincide with the probability calculated theoretically, instead the numerical value of $Fs(R_i)$ is not identical to $P(R_i)$. Why ? Is there perhaps a systematic error in the experiment or is the theory wrong ?

The relationship $R_i$ at level $n$ works by means of its subrelationships at level $(n + 1)$, however these occur randomly. They do not follow a procedure and a finite number of tests does not allow the subrelationships of $R_i$ to maintain their dynamical contribution to $R_i$. Symmetrically the subrelationships of $NOT\ R_i$ are not proportional to $P(NOT\ R_i)$. Every finite sample of tests unbalances both $R_i$ and $NOT\ R_i$. The occurrences of one group are lower to what they ought to be and the occurrences of the other are greater. The relative frequencies appear in favor of one group and in detriment of the other. $Fs(R_i)$ and $Fs(NOT\ R_i)$ are necessarily unreliable and disagree $P(R_i)$ and $P(NOT\ R_i)$ which we calculate. We conclude the correct trial of probability must be extended over the universe where the subrelationships of $R_i$ and of *NOT* $R_i$ do not undergo limitations. The ideal experimentation of $P(R_i)$, which excludes any deforming influence and provides the right $Fs(R_i)$, requires the number $Gs$ of tests be infinite

$$Gs = \infty \qquad (31)$$

In this situation the theoretical and the experimental measurement of $R_i$ coincide

$$|Fs(R_i) - P(R_i)| = 0 \qquad (32)$$

This ideal experiment is impossible to conduct and so we go as far as is possible. We define this approximation using the limit

$$\lim_{Gs \to \infty} |Fs(R_i) - P(R_i)| = 0 \qquad (33)$$

We make the limit explicit and we claim that, given a high number $N$, there is a value $Gs$

$$Gs > N \qquad (34)$$

such that

$$|Fs(R_i) - P(R_i)| < \frac{1}{Gs} \qquad (35)$$

In practice, if we carry out the tests a "sufficiently high" number of times then the difference between the frequency and the probability will be less to the "small" number $1/Gs$. The limit ensures as fine a result as we want and this is enough for correctness of the method. Equation (33) proves that the probability defined by (20), (21), (22) is verifiable in fact and confirms that the present theory has substance.

Note that the limit (33) is compatible with the *empirical law of chance* or *law of great numbers* but does not define probability. It simply asserts how probability can be verified in practice and does not give rise to the conceptual difficulties that we find in the frequentists' definition [12].

**Remark 6.6:** The limit (33) holds that the higher the number of trials the more frequency gets close to probability. Vice versa, the smaller is the sample, the less the test of probability is reliable and the maximum deviation emerges in one test. In the single event one subrelationship of the level $(n + 1)$ fires the process and this subrelationship pertains to $R_i$ or otherwise pertains to *NOT* $R_i$. In both the alternatives $Fs(R_i)$ disagrees completely from the probability, which should be decimal.

| **Gs** | 1 | > N | ∞ |
|---|---|---|---|
| **Fs** | wrong | approximate | right |

$$\qquad (36)$$

This spectrum clarifies the significance of $Fs(R_i)$ and also of $P(R_i)$; in fact any scientific measure takes its significance under the precise conditions pertaining to it. A parameter is not valid everywhere and anywhere, but it does only under precise conditions. A fair example can clarify this statement. In Mechanics we define the force $f$ as the factor causing the acceleration $a$ to the mass

$$f = m \cdot a \qquad (37)$$

This expression is true within the inertial system which is stationary or moving straight on and steadily. The mass *m* accelerates after the force in accordance with (36) in the inertial system, conversely the body can accelerate without any force in the non-inertial reference. The force cannot be tested and definition (37) has no validity when system is not inertial. The same criterion applies to probability, which takes on significance only under the experimental conditions pertaining to it, and has no objective meaning outside of this context. The experimental conditions of probability are expressed by the limit (33) and are somewhat complex. We have not two alternative and mutually exclusive reference systems, inertial and non-inertial, conversely we have the continuous spectrum (36). *Probability is correctly experimented and thus takes on a right and objective significance* if

$$Gs = \infty \quad (38)$$

This is unattainable and we use a large sample

$$Gs > N \quad (39)$$

The higher is the test number and the more objective is $P(R_i)$ vice versa it loses significance as much as *Gs* decreases. The probability is absolutely not valid if

$$Gs = 1 \quad (40)$$

Universal experience exhibits that the decimal value of the probability diverges from the frequency which is unit or otherwise zero

$$P(R_i) \neq F(R_i) \quad (41)$$

Ancient philosophers were aware of this open contradiction and along several centuries refused any indeterministic reasoning. The developments of natural sciences in the nineteenth century highlighted the relevance of the probability (see Remark 6.4) and solicited specialists to calculate $P(R_i)$ even in one trial. However this number is absolutely incorrect in point of science since it cannot be validated in the physical world. Probability can orientate only the personal expectation, namely probability takes a mere *subjective significance* [13].

| *Gs* | 1 | > N | ∞ |
|---|---|---|---|
| *Fs* | wrong | approximate | right |
| *P* | subjective | objective | |

$$(42)$$

$F(R_i)$ and $P(R_i)$ have a correct and objective meaning when they refer to the complete inductive base. As the number of experiments decrease so the precision of $F(R_i)$ decreases and the objectivity of $P(R_i)$ decreases progressively. In (40) the numerical value of $F(R_i)$ is systematically wrong and the value of $P(R_i)$ is subjective.

## VII. CONCLUSIONS

This paper relates our explicit criticism on the modeling of the probabilistic event. We highlight that Pascal's conjecture has been accepted uncritically and the argument of the probability was not studied sufficiently in depth. Non-marginal errors emerge in the relation between the formal models and the physical experiments.

We put forward the algebraic structure that presents several advantages. It displays any component occurring in the world. It suggests a definition of the probability which has a precise significance on the theoretical and experimental planes. The complex relations between $P(R_i)$ and the corresponding physical measurement $F(R_i)$ appear open. Moreover they give an answer to the everlasting debate upon the subjective and objective interpretations of probability.

This contribution is taken out of a larger study on the foundations of probability[14] which expands the analysis of physical events. It took its initial inspiration from the software system analysis with which we were familiar.